\documentclass{amsart}
\usepackage[foot]{amsaddr}
\usepackage{amsmath}

\newtheorem{theorem}                   {Theorem}

\begin{document}

\title{Bootstrap percolation in random $k$-uniform hypergraphs}

\thanks{\textsuperscript{1}supported by DFG KA 2748/3-1 and Austrian Science Fund (FWF): W1230, P26826}
\thanks{\textsuperscript{2}supported by NAWI Graz and Austrian Science Fund (FWF): W1230, P26826} 
\thanks{\textsuperscript{3}supported by DFG KA 2748/3-1 and Austrian Science Fund (FWF): P26826}
\thanks{\textrm{DOI:} https://doi.org/10.1016/j.endm.2015.06.081. Copyright \textcopyright 2015. This manuscript version is made available under the CC-BY-NC-ND 4.0 license http://creativecommons.org/licenses/by-nc-nd/4.0/}

\author[M.~Kang, C.~Koch and T.~Makai]{Mihyun Kang\textsuperscript{1}, Christoph Koch\textsuperscript{2} and Tam\'as Makai\textsuperscript{3}}

\email{\{kang,ckoch,makai\}@math.tugraz.at}
\address{Institute of Optimization and Discrete Mathematics, Graz University of Technology, 8010 Graz, Austria}

\date{}

\begin{abstract}
We investigate bootstrap percolation with infection threshold $r> 1$ on the binomial $k$-uniform random hypergraph $H_k(n,p)$ in the regime $n^{-1}\ll n^{k-2}p \ll n^{-1/r}$, when the initial set of infected vertices is chosen uniformly at random from all sets of given size. We establish a threshold such that if there are less vertices in the initial set of infected vertices, then whp only a few additional vertices become infected, while if the initial set of infected vertices exceeds the threshold then whp almost every vertex becomes infected. In addition, for $k=2$, we show that the probability of failure decreases exponentially.

\end{abstract}

\maketitle

\section{Introduction}

Bootstrap percolation on a hypergraph with infection threshold $r\geq 1$ is a deterministic infection process which evolves in rounds. In each round every vertex has exactly one of two possible states: it is either infected or uninfected. We denote the set of initially infected vertices by $\mathcal{A}_r(0)$. We say that a vertex $u$ is a neighbour of $v$ if there exists an edge containing both $u$ and $v$. In each round of the process every uninfected vertex $v$ becomes infected if it has at least $r$ infected neighbours, otherwise it remains uninfected. Once a vertex has become infected it remains infected forever. The process stops once no more vertices become infected and we denote this time step by $T$. The final infected set is denoted by $\mathcal{A}_r(T)$.

Bootstrap percolation was introduced by Chalupa, Leath, and Reich \cite{bootstrapintr} in the context of
magnetic disordered systems. Since then bootstrap percolation processes (and extensions) have been used to describe several complex phenomena: from neuronal activity \cite{MR2728841,inhbootstrap} to the dynamics of the Ising model at zero temperature \cite{Fontes02stretchedexponential}.

In the context of social networks, bootstrap percolation provides a prototype model for the spread of ideas. In this setting infected vertices represent individuals who have already adopted a new belief and a person adopts a new belief if at least $r$ of his acquaintances have already adopted it.

On the $d$-dimensional grid $[n]^d$ bootstrap percolation has been studied by 
Balogh, Bollob{\'a}s, Duminil-Copin, and Morris \cite{MR2888224}, when the initial infected set contains every vertex independently with probability $p$. 
For the size of the final infection set they showed the existence of a sharp threshold. More precisely, they established the threshold probability $p_\mathrm{c}$, such that if $p\leq (1-\varepsilon )p_\mathrm{c}$, then the probability that every vertex in $[n]^d$ becomes infected tends to 0, as $n\rightarrow\infty$, while if $p\geq (1+\varepsilon )p_\mathrm{c}$, then the probability that every vertex in $[n]^d$ becomes infected tends to one, as $n\rightarrow\infty$.

Bootstrap percolation has also been studied for several random graph models. For instance Amini and Fountoulakis \cite{bootpower} considered the Chung-Lu model \cite{MR1955514} where the vertex weights follow a power law degree distribution and the presence of an edge $\{u,v\}$ is proportional to the product of the weights of $u$ and $v$. Taking into account that in this model a linear fraction of the vertices have degree less than $r$ and thus at most a linear fraction of the vertices can become infected the authors proved the size of the final infected set $\mathcal{A}_r(T)$ exhibits a phase transition. 


Janson,  \L uczak, Turova, and Vallier \cite{MR3025687} analysed bootstrap percolation on the binomial random graph $G(n,p)$ where every edge appears independently with probability $p$. For $r\geq 1$ and $n^{-1}\ll p \ll n^{-1/r}$ they showed that there is a threshold such that if the initial number of infected vertices is below the threshold, then the process infects only a few additional vertices and if the initial number of infected vertices exceeds the threshold, then almost every vertex becomes infected.

In this paper we investigate the binomial random hypergraph $H_k(n,p)$, where every edge ($k$-tuple of vertices) is present independently with probability $p$. We choose the initial infected set uniformly at random and consider bootstrap percolation with infection threshold $r> 1$ in the regime $n^{-1}\ll n^{k-2}p \ll n^{-1/r}$.
%
The main contribution of this paper are:
\begin{itemize}
\item strengthening of the result in \cite{MR3025687}, by showing that the failure probability decreases {\em exponentially} (Theorem~\ref{thm:graph});
\item extension of the original results from graphs to hypergraphs (Theorem~\ref{thm:hypergraph}).
\end{itemize}

\section{Main Results}

We extend the following result, which was originally proved in \cite{MR3025687}, to $H_k(n,p)$: Consider bootstrap percolation with infection threshold $r$ on $G(n,p)$, where $n^{-1}\ll p \ll n^{-1/r}$. There is a threshold $b_r=b_r(n,p)$ such that if $|\mathcal{A}_r(0)|\le(1-\varepsilon)b_r$, then with probability tending to one as $n\rightarrow \infty$ (whp for short)  only a few additional vertices become infected, while if $|\mathcal{A}_r(0)|\ge(1+\varepsilon)b_r$, then whp almost every vertex in the process becomes infected. 
%
For integers $k\geq 2$ and $r> 1$ set
\[
b_{k,r}:=b_{k,r}(n,p)=\left\{ \begin{array}{ll}
\left(1-\frac{1}{r}\right)\left(\frac{(r-1)!}{n\left(\binom{n}{k-2}p\right)^r}\right)^{1/(r-1)} & \mbox{if } r > 2 \\
\frac{1}{2(2k-3)}\frac{1}{n\left(\binom{n}{k-2}p\right)^2} & \mbox{if } r = 2,
\end{array}
\right.
\]
and note that the only difference for the $r=2$ case is a $1/(2k-3)$ multiplier. Since $2k-3=1$ when $k=2$ this is consistent with the threshold in the graph case i.e.\ $b_{2,r}=b_r$. 

\begin{theorem}\label{thm:hypergraph}
For $k\geq 2$ consider bootstrap percolation with infection threshold $r> 1$ on $H_{k}(n,p)$ when $n^{-1}\ll n^{k-2} p \ll n^{-1/r}$. Assume the initial infection set is chosen uniformly at random from all sets of vertices of size $a=a(n)$. Then for any fixed $\varepsilon>0$ we have that
\begin{itemize}
\item if $a\leq(1-\varepsilon)b_{k,r}$ then whp $|\mathcal{A}_r(T)|= O(b_{k,r})$;
\item if $a\geq (1+\varepsilon)b_{k,r}$ then whp $|\mathcal{A}_r(T)|=(1+o(1))n$.
\end{itemize}
\end{theorem}

Using the methods developed for this result we also obtain a strengthened form of the result for $G(n,p)$ establishing exponentially small bounds on the failure probability.

\begin{theorem}\label{thm:graph}
Consider bootstrap percolation with infection threshold $r > 1$ on $G(n,p)$ when $n^{-1}\ll p \ll n^{-1/r}$. Assume the initial infection set is chosen uniformly at random from the set of vertices of size $a=a(n)$. Then for any fixed $\varepsilon>0$ the following holds with probability $1-\exp(-\Omega(b_{2,r}))$:
\begin{itemize}
\item if $a\leq(1-\varepsilon)b_{2,r}$, then $|\mathcal{A}_r(T)|=O(b_{2,r})$;
\item if $a\geq (1+\varepsilon)b_{2,r}$, then $|\mathcal{A}_r(T)|=(1+o(1))n$.
\end{itemize}
\end{theorem}

The proofs rely on surprisingly simple methods. When the number of vertices infected in the individual rounds is large, we apply Chebyshev's or Chernoff's inequality. However when the process dies out, these changes can become arbitrarily small. In this case we couple the infection process with a subcritical branching process which dies out very quickly.



\section{Proof outlines}

We first show the outline for the proof of Theorem~\ref{thm:hypergraph}. For brevity we will only describe the $r>2$ case in detail and comment on the differences for $r=2$ at the end. 

Start with a given set of initially infected vertices $\mathcal{A}_r(0)$ and consider the infection process round by round. At the end of round $t\geq 1$ we partition the set of vertices into $\mathcal{A}_0(t),\mathcal{A}_1(t),...,\mathcal{A}_r(t)$ where the set $\mathcal{A}_i(t)$ consists of all the vertices which have exactly $i$ infected neighbours (these are vertices in $\mathcal{A}_r(t-1)$), for $i<r$, and $\mathcal{A}_r(t)$ consists of all the vertices which have at least $r$ infected neighbours.

For every $0\leq i \leq r$ we aim to define a sequence $\{a_i(t)\}_{t\geq 0}$ in such a way that $|\mathcal{A}_i(t)|\approx a_i(t)$. We use the following initial values for the sequences $a_0(0)=n$, $a_r(0)=a$, and $a_i(0)=0$ for $0< i < r$.

As long as $|\mathcal{A}_r(t)|=o\left(\left(n^{k-2}p\right)^{-1}\right)$ the expected number of infected neighbours of a vertex is $o(1)$ and thus the typical vertex that becomes infected in this round has exactly $r$ infected neighbours in $r$ different edges. 
In round $t+1$  we determine for any uninfected vertex $v\in V\backslash \mathcal{A}_r(t)$ whether it changes its partition class. Note that this only happens if it has at least one neighbour which became infected in round $t$. 
Therefore, for $0< i \leq r$, the expected change $|\mathcal{A}_i(t+1)\backslash \mathcal{A}_i(t)|$ of the size of the partition class can be approximated by
\begin{equation}\label{eq:seqchange}
a_{i}(t+1)-a_i(t) \approx \sum_{j=1}^{i} \left(\frac{(a_r(t)-a_{r}(t-1))^j}{j!}a_{i-j}(t)\right)\left(\binom{n}{k-2}p\right)^{i-j}, 
\end{equation}
ignoring any negative terms (which correspond to vertices leaving their partition class). Similarly we assume that the number of vertices without infected neighbours does not change significantly, i.e.\ $a_0(t+1)=a_0(t)=n$.
From \eqref{eq:seqchange} we deduce
\begin{equation}\label{eq:seq}
a_i(t+1)\approx \frac{a_r(t)^i}{i!}n\left(\binom{n}{k-2}p\right)^i+a_i(0), 
\end{equation}
for $0< i \leq r$. The behaviour of the sequence depends on the size $a=|\mathcal{A}_r(0)|$ of the initial set. We will show the following: in the subcritical regime, characterised by $a\leq (1-\varepsilon)b_{k,r}$, the sequence $\{a_r(t)\}_{t\geq 0}$ converges to a value $a^*=O(b_{k,r})$ as $t\rightarrow \infty$. On the other hand, in the supercritical regime, $a\geq (1+\varepsilon)b_{k,r}$, the sequence $\{a_r(t)\}_{t\geq 0}$ tends to infinity as $t\rightarrow \infty$. 

First consider the subcritical regime. Since in this case $a_r(t)$ converges we have that the differences $\Delta(t):=a_{r}(t+1)-a_r(t)$ form a decreasing function in $t$ and show that, for any fixed $\eta>0$, there exists a $\tau$, which does not depend on $n$, such that $\Delta(\tau)\leq \eta b_{k,r}$. 
The fact that $|\mathcal{A}_i(t)|$ is concentrated around $a_i(t)$ for $t<\tau$ follows from Chebyshev's inequality. 
%

Since we are in the subcritical regime the size of the individual generations will become small and the concentration will fail. In order to avoid this we attempt to analyse the remaining steps together. Consider the forest where every vertex in $\mathcal{A}_r(\tau+1)\backslash \mathcal{A}_r(\tau)$ is a root.
Recall that in order for a vertex to become infected in round $t+1$ it must have a neighbour that got infected in round $t$. The children of a vertex $v \in \mathcal{A}_r(t+1)\backslash \mathcal{A}_r(t)$ will be the vertices $u \in \mathcal{A}_r(t+2)\backslash \mathcal{A}_r(t+1)$ which lie in an edge containing $v$ and should this relation not be unique for some vertex $u$, $u$ is assigned arbitrarily to one of the candidates. Clearly every vertex of $A_r(T)\setminus A_r(\tau)$ is contained in the forest and thus the size of this forest matches the number vertices which got infected after round $\tau$.

Note that for every $\delta>0$ there exists a $t_0$ such that $|\mathcal{A}_i(t)|\leq (1+\delta)a_i(\tau)$, for every $0\leq i\leq r$ and $\tau<t \leq t_0$. Also up until time $t_0+1$ we have an upper coupling by a Galton-Watson branching process with $|\mathcal{A}_r(\tau+1)\backslash \mathcal{A}_r(\tau)|$ roots and offspring distribution $\sum_{j=0}^{r-1}\mathrm{Bin}((1+\delta)a_{r-j}(\tau),q_j)$ where
\[q_j=\binom{n}{k-2}p \frac{(\delta a_r(\tau))^{j-1}}{(j-1)!}\left(\binom{n}{k-2}p\right)^{r-j-1}.\]
For small enough $\delta$ the expected number of offspring in one step is \linebreak $\sum_{j=0}^{r-1}(1+\delta)a_{r-j}(\tau) q_j<1$ and therefore this is a subcritical process, i.e.\ it dies out with probability 1.
For every $t$ we have that $|\mathcal{A}_r(t)|\leq (1+\delta)a_r(\tau)$ implies $|\mathcal{A}_i(t+1)|\leq (1+\delta)a_i(\tau)$, for all $0\leq i < r$, by \eqref{eq:seq}, and thus it is enough to show that $|\mathcal{A}_r(T)|\leq (1+\delta)a_r(\tau)$. Due to the upper coupling with the branching process we have that the probability that $|\mathcal{A}_r(T)|>(1+\delta)a_r(\tau)$ is dominated by the probability that the total size of the branching process exceeds $\delta a_r(\tau)$. However for properly chosen $\eta,\delta>0$ the probability that the total size of the branching process exceeds $\delta a_r(\tau)$ is sufficiently small. Therefore we have that there are at most $(1+\delta)a_{r}(\tau)$ infected vertices in total.


Now for the supercritical case. Recall that \eqref{eq:seqchange} and \eqref{eq:seq} hold when \linebreak[4] $a_r=o\left(\left(n^{k-2}p\right)^{-1}\right)$. Again we consider the differences $\Delta(t)=a_{r}(t+1)-a_r(t)$. Although at the beginning of the process the values of $\Delta(t)$ decrease there exists a value $t_1$ not depending on $n$ such that for every $t>t_1$ we have that $\Delta(t+1)>\Delta(t)$. In fact there exists a $t_2$ not depending on $n$ such that for $t\geq t_2$ we have that $\Delta(t+1)>2\Delta(t)$. 
Therefore the probability of non-concentration is dominated by a geometric sequence and applying the union bound gives us concentration as long as $a_r(t)=o\left(\left({n}^{k-2}p\right)^{-1}\right)$. 
When $a_r(t)=\Omega\left(\left({n}^{k-2}p\right)^{-1}\right)$ the expected number of neighbours is $\Omega(1)$ and thus our approximation in \eqref{eq:seqchange} does not hold any more. Refining these approximations shows that at most 2 rounds are required for almost every vertex to become infected, with $\Theta(n)$ vertices becoming infected in every required step.

Recall that for $r>2$ the typical vertex became infected when it was contained in $r$ different edges each containing a different infected vertex. When  $r=2$ this is equivalent to finding two intersecting edges each containing a different infected vertex. However unlike the $r>2$ case finding two such edges in step $t$ implies that every vertex in these edges is infected by step $t+1$. Two intersecting edges typically overlap in exactly one vertex and thus finding such an edge pair implies that $2k-3$ vertices will become infected, not just one. Taking this into account gives us the modified bound on the threshold.

The proof of Theorem~\ref{thm:graph} is analogous. In the random graph case, in round $t$ of the process only those edges are examined which contain exactly one vertex from $\mathcal{A}(t)\backslash \mathcal{A}(t-1)$ and no vertices from $\mathcal{A}(t-1)$. Since each of these edges can contain at most one uninfected vertex the behaviour of the individual vertices is independent. Thus we can replace Chebyshev's inequality with Chernoff's inequality and achieve a stronger bound on the failure probability.

\bibliographystyle{plain}
\bibliography{referencesEuroComb}

\end{document}